\input amstex\documentstyle{amsppt}  
\pagewidth{12.5cm}\pageheight{19cm}\magnification\magstep1  
\topmatter
\title From conjugacy classes in the Weyl group to representations\endtitle
\author George Lusztig and Zhiwei Yun\endauthor
\address{Department of Mathematics, M.I.T., Cambridge, MA 02139}\endaddress
\thanks{G.L. supported by NSF grant DMS-1855773; Z.Y. supported by a Packard Fellowship.}\endthanks
\endtopmatter   
\document

\define\Irr{\text{\rm Irr}}

\define\mpb{\medpagebreak}

\define\lb{\linebreak}

\define\op{\oplus}
   
\define\part{\partial}
\define\emp{\emptyset}

\define\m{\mapsto}
\define\do{\dots}

\define\sub{\subset}    

\define\T{\times}

\define\nl{\newline}
\redefine\i{^{-1}}

\define\un{\underline}

\define\ot{\otimes}
\define\bbq{\bar{\QQ}_l}

\define\Hom{\text{\rm Hom}}

\define\ind{\text{\rm ind}}

\define\sgn{\text{\rm sgn}}
\define\tr{\text{\rm tr}}

\define\g{\gamma}

\define\e{\epsilon}

\redefine\o{\omega}
\define\p{\pi}

\define\r{\rho}
\define\s{\sigma}

\define\th{\theta}

\redefine\l{\lambda}
\define\z{\zeta}
\define\x{\xi}

\define\Th{\Theta}

\define\Ph{\Phi}
\define\Ps{\Psi}

\define\boc{\bold c}

\define\qq{\bold q}

\define\FF{\bold F}

\define\NN{\bold N}

\define\QQ{\bold Q}

\define\ZZ{\bold Z}

\define\cb{\Cal B}

\define\ce{\Cal E}
\define\cf{\Cal F}

\define\ch{\Cal H}

\define\car{\Cal R}

\define\cu{\Cal U}

\define\sha{\sharp}

\head 1. Definition of the map $\Ps$\endhead
\subhead 1.1\endsubhead
Let $G$ be a connected reductive algebraic group defined and split over the finite field $\FF_p$ with $p$
elements where $p$ is a prime number which we assume to be good for $G$. 
Let $W$ be the Weyl group of $G$
(a Coxeter group with length function $w\m|w|$),
let $\un W$ be the set of conjugacy classes of $W$ and let $\car_W$ be the Grothendieck group of (finite 
dimensional) representations of $W$ over $\bbq$; here $l$ is a fixed prime number, $l\ne p$.
Let $\cu$ be the variety of unipotent elements of $G$. Let $\un\cu$ be the set of $G$-conjugacy classes of
unipotent elements. Let $\cb$ be the variety of Borel subgroups of $G$. 
We fix $B^+\in\cb$ such that $B^+$ is defined over $\FF_p$.
In \cite{L11} a (surjective) map $\Ph:\un W@>>>\un\cu$ was defined based on the study of intersections of
unipotent classes in $G$ with $(B^+,B^+)$ double cosets in $G$.
Let $\qq$ be an indeterminate. In this paper we define a map $\Ps:\un W@>>>\ZZ[\qq]\ot\car_W$ whose image is
denoted by $Z$. We also define a (surjective) map $\Th:Z@>>>\un\cu$ such that $\Ph(C)=\Th(\Ps(C))$ for any
$C\in\un W$. Thus, $\Ps$ is a refinement of $\Ph$. (The definition of $\Ps$ again involves 
the study of intersections of unipotent classes in $G$ with $(B^+,B^+)$ double cosets in $G$.)
We also describe $\Ps$ explicitly for $G$ of low rank.

\subhead 1.2\endsubhead
For $(B,B')\in\cb\T\cb$ we denote by $pos(B,B')\in W$ the
relative position of $B,B'$ (see \cite{DL}).
For $w\in W$ let $B^+wB^+=\{g\in G;pos(B^+,gB^+g\i)=w\}$. Let $F_1:G@>>>G$ be the Frobenius map relative to the
$\FF_p$-structure.
Let $q=p^e$ ($e\ge1$) and let $F=F_1^e:G@>>>G$. This induces Frobenius maps 
$\cu@>>>\cu$, $\cb@>>>\cb$ denoted again 
by $F$. For $y\in W$, we set $X_y=\{B\in\cb;pos(B,F(B))=y\}$ (see \cite{DL}). If $i\in\ZZ$, the 
$\bbq$-cohomology space with compact support $H^i_c(X_y)$ has a natural action of the finite group $G^F$ (see 
\cite{DL}); here $?^F$ denotes the fixed point set of $F:?@>>>?$. Let $R^1_y=\sum_i(-1)^iH^i_c(X_y)$,
a virtual representation of $G^F$. For $g\in G^F$, $\tr(g,R^1_y)$ is an integer, see \cite{DL, 3.3}. For 
$w\in W,y\in W$ we set
$$a_{w,y}=\sum_{u\in\cu^F\cap(B^+wB^+)}\tr(u,R^1_y)\in\ZZ.\tag a$$

\proclaim{Proposition 1.3} Let $w\in W$. There is a unique element $\z_{w,q}\in\car_W$ such that for any 
$y\in W$ we have $\tr(y,\z_{w,q})=a_{w,y}$.
\endproclaim
For $u\in\cu$ we set $\cb_u=\{B\in\cb_u\in B\}$. It is known (Springer) that $W$ acts naturally on the
$l$-adic cohomology space $H^i(\cb_u)$. 
Let $\Irr(W)$ be a set of representatives for the isomorphism classes of irreducible representations of $W$
over $\bbq$. We can identify $H^i(\cb_u)=\op_{E\in\Irr(W)}E\ot H^i(\cb_u)_E$
where $H^i(\cb_u)_E=\Hom_W(E,H^i(\cb_u))$.
From \cite{L90} (or \cite{K} when $p\gg0$), for $u\in\cu^F$ and $y\in W$ we have
$$\align&\tr(u,R^1_y)\\&=
\sum_i(-1)^i\tr(yF,H^i(\cb_u))=\sum_{E\in\Irr(W)}\tr(y,E)\tr(F,H^*(\cb_u)_E).\tag a\endalign$$
 Here $F:H^i(\cb_u)@>>>H^i(\cb_u)$ and $F:H^i(\cb_u)_E@>>>H^i(\cb_u)_E$ are induced by
the restriction of $F:\cb@>>>\cb$ to $\cb_u$ and $\tr(F,H^*(\cb_u)_E)$ is defined to be
$$\sum_i(-1)^i\tr(F,H^i(\cb_u)_E).$$ Thus we have
$$a_{w,y}=\sum_{u\in\cu^F\cap(B^+wB^+)}\sum_{E\in\Irr(W)}\tr(y,E)\tr(F,H^*(\cb_u)_E).$$
For $E\in\Irr(W)$ we set 
$$\z_{w,q;E}=\sum_{u\in\cu^F\cap(B^+wB^+)}\tr(F,H^*(\cb_u)_E)\tag b$$
so that
$$a_{w,y}=\sum_{E\in\Irr(W)}\z_{w,q;E}\tr(y,E)$$
for any $y\in W$. Since $a_{w,y}\in\ZZ$ and the matrix $(\tr(y,E))$ (with $y$ running through representatives
for the conjugacy classes in $W$, and $E\in\Irr(W)$) has integer entries and nonzero determinant, we see that
$\z_{w,q;E}\in\QQ$ for any $E\in\Irr(W)$. On the other hand from the definition we see see that $\z_{w,q;E}$ is an
algebraic integer. It follows that 
$$\z_{w,q;E}\in\ZZ.\tag c$$
We set $\z_{w,q}=\sum_{E\in\Irr(W)}\z_{w,q;E}E\in\car_W$. This proves the existence statement in the proposition.
The uniqueness is obvious.

\subhead 1.4\endsubhead
Let $\cf$ be the vector space of functions $\cb^F@>>>\bbq$. For $w\in W$ we define a linear map
$T_w:\cf@>>>\cf$ by $T_w(f)(B)=\sum_{B'\in\cb^F;pos(B,B')=w}f(B')$. Let $\ch$ be the subspace of $End(\cf)$
with basis $\{T_w;w\in W\}$; this is a subalgebra of $End(\cf)$ (the Hecke algebra). Now $G^F$ acts on $\cf$
by $g:f\m f'$ where $f'(B)=f(g\i Bg)$. This action commutes with the $\ch$-action so that we can identify
$\cf=\op_{\ce\in\Irr(\ch)}\ce\ot\cf_\ce$ where $\cf_\ce=\Hom_\ch(\ce,\cf)$ and $\Irr(\ch)$ is a set of 
representatives for the simple $\ch$-modules; here $\cf_\ce$ is an irreducible $G^F$-module.
For $u\in\cu^F$ we have  
$$\sha(B\in\cb^F;pos(B,uBu\i)=w)=\tr(uT_w,\cf)=\sum_{\ce\in\Irr(\ch)}\tr(T_w,\ce)\tr(u,\cf_\ce)$$
hence for $y\in W$ we have
$$\align&\sha(\cb^F)a_{w,y}=\sum_{u\in\cu^F,B\in\cb^F;pos(B,uBu\i)=w}\tr(u,R^1_y)\\&
=\sum_{\ce\in\Irr(\ch)}\sum_{u\in\cu^F}\tr(u,\cf_\ce)\tr(u,R^1_y)\tr(T_w,\ce).\tag a\endalign$$
From (a) and 1.3(a) for any $w\in W,E\in\Irr(W)$ we have:
$$\sha(\cb^F)\z_{w,q;E}
=\sum_{u\in\cu^F}\sum_{\ce\in\Irr(\ch)}\tr(u,\cf_\ce)\tr(F,H^*(\cb_u)_E)\tr(T_w,\ce).\tag b$$

\subhead 1.5\endsubhead
Let $w\in W,E'\in\Irr(W)$. We set
$$c_{w,E',q}=\sha(W)\i\sum_{\ce\in\Irr(\ch)}\sum_{z\in W}(\cf_\ce:R^1_z)_{G^F}\tr(z,E')\tr(T_w,\ce)
\in\QQ.\tag a$$ 
Here $(:)_{G^F}$ is the standard inner product of virtual representations of $G^F$.
We shall now regard $q$ as variable. We show:

(b) {\it $c_{w,E',q}$ is the value at $q$ of a polynomial
$\boc_{w,E'}(\qq)$ in $\qq$ with rational coefficients independent of $q$.}
\nl
We can assume that $G$ is adjoint simple. The only part of the right hand side of (a) which depends on $q$ is
$\tr(T_w,\ce)$. This is known to be the value at $q$ of a polynomial in $\qq$ with integral coefficients
independent of $q$ except when $G$ is of type $E_7$, $\dim E'=512$ or $G$ is of type $E_8$, $\dim E'=4096$.
Assume now that $G$ is of type $E_7$, $\dim E'=512$. Let $\ce',\ce''$ be the two objects of $\Irr(\ce)$ which
have dimension 512. We have $c_{w,E',q}=(1/2)(\tr(T_w,\ce')+\tr(T_w,\ce''))$ (see the proof of \cite{L84, 7.12}).
Although the quantities $\tr(T_w,\ce'),\tr(T_w,\ce'')$ are separately not necessarily polynomials in $q$
their sum is. This proves (b) in the case where $G$ is of type $E_7,\dim E'=512$. The proof in the case
where $G$ is of type $E_8, \dim E'=4096$, is entirely similar.

\subhead 1.6\endsubhead
For $\ce\in\Irr(\ch)$ there is a well defined class function $\x_\ce$ on $G^F$ such that
$$\tr(g,\cf_\ce)=\sha(W)\i\sum_{z\in W}(\cf_\ce:R^1_z)_{G^F}\tr(g,R^1_z)+\x_\ce(g)\tag a$$
for all $g\in G^F$ and 

(b) $\x_\ce$ is orthogonal to the character of any virtual representation $R_T^\th$ of \cite{DL}.
\nl
From (b) we deduce using results in \cite{L86} that for any $E\in\Irr(W)$ we have 
$$\sum_{u\in\cu^F}\x_\ce(u)\tr(F,H^*(\cb_u)_E)=0.$$
Using (a) and 1.4(b) we see that for any $w\in W,E\in\Irr(W)$ we have
$$\align&\sha(\cb^F)\z_{w,q;E}\\&
=\sha(W)\i\sum_{u\in\cu^F}\sum_{\ce\in\Irr(\ch)}\sum_{z\in W}(\cf_\ce:R^1_z)_{G^F}\tr(u,R^1_z)
\tr(F,H^*(\cb_u)_E)\tr(T_w,\ce).\endalign$$
Using 1.3(a) we deduce
$$\sha(\cb^F)\z_{w,q;E}=\sum_{E'\in\Irr(W)}c_{w,E',q}\sum_{u\in\cu^F}
\tr(F,H^*(\cb_u)_{E'})\tr(F,H^*(\cb_u)_E).\tag c$$
The right hand side of (c) is the value at $\qq=q$ of a polynomial $\p(\qq)$ (with 
the coefficients of $\p$ being rational numbers independent of $q$).
(This follows from 1.5(b) and the known properties of $\tr(F,H^i(\cb_u)_E)$, see \cite{DLP}.)

\subhead 1.7\endsubhead
Let $w\in W,E\in\Irr(W)$. It is known that 

(a) {\it $H^i(\cb_u)_E$ can be interpreted as the stalk at $u\in\cu$ of
an intersection cohomology complex $K_E$ on the closure of a unipotent class of $G$.}
\nl
(This was stated in \cite{L81, \S3, Conj.2} and proved in \cite{BM}.)
We use this and the Grothendieck trace formula to rewrite the definition of $\z_{w,q;E}$ in the form
$$\z_{w,q;E}=\sum_i(-1)^i\tr(F,H^i(X,K_E|_X))\tag b$$
where $X=\cu\cap(B^+wB^+)$ and the map $H^i(X,K_E|_X)@>>>H^i(X,K_E|_X)$ induced by $F:G@>>>G$ is denoted again 
by $F$. Replacing $q$ by $p^s$ with $s=1,2,\do$ we see that
$$\z_{w,p^s;E}=\sum_{\l\in V}n_\l\l^s$$
where $V$ is a finite subset of $\bbq-\{0\}$ and $n_\l$ are nonzero integers independent of $s$. 
It follows that
$$(\sum_{v\in W}p^{|v|s})\z_{w,p^s;E}=\sum_{\l'\in V'}n'_{\l'}\l'{}^s$$
where $V'$ is a finite subset of $\bbq-\{0\}$ and $n'_{\l'}$ are nonzero integers independent of $s$. 
By the results in 1.6 we have also
$$(\sum_{v\in W}p^{|v|s})\z_{w,p^s;E}=\sum_{m\in[0,N]}n''_mp^{sm}$$
for some $N\ge1$ where $n''_m\in\QQ$ are independent of $s$. We deduce that
$$\sum_{\l'\in V'}n'_{\l'}\l'{}^s=\sum_{m\in[0,N]}n''_mp^{sm}$$
for $s=1,2,\do$. This forces $V'$ to be a subset of $\{1,p,p^2,\do\}$. Thus we have the following result.

(b) {\it There exists a polynomial $\p(\qq)$ in $\qq$ with integer coefficients independent of $s$
such that 

$(\sum_{v\in W}p^{|v|s})\z_{w,p^s;E}=\p(p^s)$ 

for $s=1,2,\do$.}

\subhead 1.8\endsubhead
For $C\in\un W$ we denote by $C_{min}$ the set of elements of minimal length in $C$.
Let $\r$ be the reflection representation of $W$. For $w\in W$ let $n_w=\det(1-w,\r)$. We have $n_w\ge0$.
We say that $w$ is elliptic if $n_w>0$.
Let $\un W_{el}$ be the set of all $C\in\un W$ such that for some/any $w\in C$, $w$ is elliptic.

Let $Z_G$ be the center of $G$. Let $\nu=\dim\cb$. Let $r$ be the rank of $G/Z_G$.

In this subsection we fix $C\in\un W_{el}$ and $w\in C_{min}$. Let 
$n'_w\ge1$ be the part prime to $p$ of $n_w$. According to \cite{L11, 5.2}, for any $g\in B^+wB^+$,

(a) {\it  the group $\{b\in B^+;bgb\i=g\}/Z_G$ is finite abelian of order dividing $n'_w$.}
\nl
(In the first line of \cite{L11, 5.1} one should add the sentence: {\it We fix $w\in C_{min}$}.) We show:

(b) {\it For any $E\in\Irr(W)$ we have $q^{-\nu}(q-1)^{-r}n'_w\z_{w,q;E}\in\ZZ$.}
\nl
We have $\z_{w,q;E}=\sum_S\sha(S)\tr(F,H^*(\cb_{u_S})_E)$ where $S$ runs over the set of orbits of 
of $B^{+F}$ acting on $\cu^F\cap B^+wB^+$ by conjugation and for each such $S$, $u_S$ is an element of $S$.
Here for any $S$ in the sum, $\tr(F,H^*(\cb_{u_S})_E)$ is an algebraic integer and, by (a),
$\sha(S)\in q^\nu(q-1)^rn'_w{}\i\ZZ$. Hence $q^{-\nu}(q-1)^{-r}n'_w\z_{w,q;E}$ is an algebraic integer. By 1.3(c),
this is also a rational number, hence an integer. This proves (b).

\proclaim{Proposition 1.9}Let $C\in\un W_{el}$, $w\in C_{min}$. Let $E\in\Irr(W)$.
Then \lb $\z_{w,q;E}/(q^\nu(q-1)^r)\in\ZZ$.
\endproclaim
Let $\p_0(\qq)=(\sum_{v\in W}\qq^{|v|s})\qq^\nu(\qq-1)^r$. This is a monic polynomial in $\qq$ of degree
$2\nu+r$ with integer coefficients. Let $\p(\qq)$ be as in 1.7(b). We have
$\p(\qq)=\p_1(\qq)\p_0(\qq)+\p_2(\qq)$ where $\p_1(\qq),\p_2(\qq)$ are 
polynomials with integer coefficients and $\p_2(\qq)$ is either $0$ or has degree $<2\nu+r$. By 1.7(b) and
1.8(b) for $s=1,2,\do$ we have $\p(q^s)n'_w/\p_0(w^s)\in\ZZ$ 
hence $\p_1(q^s)n'_w+\p_2(q^s)n'_w/\p_0(q^s)\in\ZZ$ so that 
$\p_2(q^s)n'_w/\p_0(q^s)\in\ZZ$. If $\p_2(\qq)\ne0$ this is impossible for large $s$ since $\deg\p_2<\deg \p_0$.
We see that $\p(\qq)=\p_1(\qq)\p_0(\qq)$.
Setting $\qq=q$ we see that $\z_{w,q;E}/(q^\nu)(q-1)^r)=\p_1(q)\in\ZZ$. The proposition is proved.

This proof shows also that for $s=1,2,\do$, 

(a) $\z_{w,p^s;E}/(p^{s\nu})(p^s-1)^r)=\p_1(p^s)$ where $\p_1(\qq)$ is a polynomial in $\qq$ with
integer coefficients independent of $s$. 

\subhead 1.10\endsubhead
Let $P$ be a parabolic subgroup of $G$ containing $B^+$. Let $L$ be an $F$-stable Levi subgroup of $P$
and let $U$ be the unipotent radical of $P$. Then $B^+_L:=B^+\cap L$ is an $F$-stable Borel subgroup of $L$
and we have $B^+=B^+_LU$. Let $\cu_L$ be the set of unipotent elements of $L$.
Let $W_L$ be the Weyl group of $L$. We can view $W_L$ as a subgroup of $W$ and as an indexing set for the 
$(B^+_L,B^+_L)$ double cosets of $L$, so that for $w\in W_L$ the double coset $B^+_LwB^+_L$
satisfies $(B^+_LwB^+_L)U=B^+wB^+$.
For $y'\in W_L$ we define $X_{y',L}$ in the same way as $X_y$, but replacing $G,y$ by $L,y'$; then
$R^1_{y',L}=\sum_i(-1)^iH^i_c(X_{y',L}))$ is naturally a (virtual) $L^F$-module.

Let $\s$ be an irreducible $L^F$-module. We can view $\s$ as a $P^F$-module on which $U^F$ acts trivially. Let $y\in W$.
The following identity is a reformulation of a special case of a result in \cite{DL1}:
$$(\ind_{P^F}^{G^F}(\s):R^1_y)_{G^F}=\sha(W_L)\i\sum_{z\in W;zyz\i\in W_L}(\s:R^1_{zyz\i,L})_{L^F}.\tag a$$
Here $(:)_{L^F}$ denote the standard inner product of virtual representations of $L^F$.
The left hand side of (a) is equal to $(\s:\sum_i(-1)^iH^i_c(X_y)^{U^F})_{L^F}$
where $H^i_c(X_y)^{U^F}$ is the space of $U^F$-invariants on $H^i_c(X_y)$ viewed as an $L^F$-module.
It follows that 
$$\sha(W_L)\sum_i(-1)^iH^i_c(X_y)^{U^F}=\sum_{z\in W;zyz\i\in W_L}R^1_{zyz\i,L}$$
as virtual $L^F$-modules.

For $w\in W_L,y'\in W_L$ we define $a_{w,y';L}$ as in 1.2(a), in terms of $L$ instead of $G$ that is,
$$a_{w,y';L}=\sum_{u\in\cu_L^F\cap(B^+_LwB^+_L)}\tr(u,R^1_{y',L})\in\ZZ.$$
Now let $w\in W_L,y\in W$. We have
$$\align&a_{w,y}=\sum_{u\in\cu^F\cap(B^+_LwB^+_L)U}\tr(u,R^1_y)\\&
=\sum_{(u',u'')\in\cu^F_L\T U^F; u'\in B^+_LwB^+_L}\tr(u'u'',R^1_y)\\&
=\sha(U^F)\sum_{u'\in\cu_L^F(B^+_LwB^+_L)}\sum_i(-1)^i\tr(u',H^i_c(X_y)^{U^F})\\&
=\sha(W_L)\i\sha(U^F)\sum_{u'\in\cu_L^F(B^+_LwB^+_L)}\sum_{z\in W;zyz\i\in W_L}
\tr(u',R^1_{zyz\i,L})\\&=\sha(W_L)\i\sha(U^F)\sum_{z\in W;zyz\i\in W_L}a_{w,zyz\i;L}.\tag b\endalign$$
Let $\car_{W_L}$ be the Grothendieck group of $W_L$-modules over $\bbq$. We define $\z_{w,L,q}\in\car_{W_L}$ as 
in 1.3 in terms of $L$
instead of $G$. Using (b) we have 
$$\tr(y,\z_{w,q})=\sha(W_L)\i\sha(U^F)\sum_{z\in W;zyz\i\in W_L}\tr(zyz\i,\z_{w,L,q})$$
that is,
$$\z_{w,q}=\sha(U^F)\ind_{W_L}^W(\z_{w,L,q}).\tag c$$

\subhead 1.11\endsubhead
Let $C\in\un W$. For $w\in C_{min},w'\in C_{min}$ we show:
$$\z_{w,q}=\z_{w',q}.\tag a$$
It is enough to show that for any $y\in W$ we have $a_{w,y}=a_{w',y}$.
We write 1.4(a) for $w,y$ and for $w',y$. We see that it is enough to show that
$\tr(T_w,\ce)=\tr(T_{w'},\ce)$ for any $\ce\in\Irr(\ch)$; this follows from results in \cite{GP}.

\subhead 1.12\endsubhead
Let $C\in\un W$. Let $m(C)$ be the multiplicity of the eigenvalue $1$ of $w$ in $\r$ (see 1.8) for some/any 
$w\in C$. We define an element $\un\z_{C,q}\in\QQ\ot\car_W$ by
$$\un\z_{C,q}=\z_{w,q}/(q^\nu(q-1)^{r-m(C)})$$
where $w\in C_{min}$. This is independent of the choice of $w$ by 1.11. We show that 
$$\un\z_{C,q}\in\car_W.\tag a$$
If $m(C)=0$ (so that $C\in\un W_{el}$) this follows from 1.9.
Assume now that $m(C)>0$. In this case $C_{min}$ contains an element $w$ which is contained and elliptic in 
$W_L$ where $P,U,L,W_L$ are as in 1.10 and $P\ne G$. By 1.9 for $L$ instead of $G$ we have
$\z_{w,L,q}/(q^{\nu-\o}(q-1)^{r-m(C)})\in\car_{W_L}$ where $\sha(U)=q^\o$.
It follows that
$\ind_{W_L}^W(\z_{w,L,q}/(q^{\nu-\o}(q-1)^{r-m(C)}))\in\car_W$. Using 1.10(c) we deduce
$\z_{w,q}/(q^\o q^{\nu-\o}(q-1)^{r-m(C)})\in\car_W$. This proves (a).

\subhead 1.13\endsubhead
Let $C\in\un W$. From 1.9(a) and the definitions we see that there is a well defined element
$\Ps(C)\in\ZZ[\qq]\ot\car_W$  such that for any $s=1,2,\do$, the specialization $\Ps(C)|_{\qq=p^s}\in\car_W$
is equal to $\un\z_{C,p^s}$ (as in 1.12). We can write $\Ps(C)=\sum_{i\ge0}\qq^i\Ps_i(C)$
where $\Ps_i(C)\in\car_W$ are zero for $i\gg0$. Thus we have a map $\Ps:\un W@>>>\ZZ[\qq]\ot\car_W$
and maps $\Ps_i:\un W@>>>\car_W$.

\head 2. Properties of the map $\Ps$\endhead
\subhead 2.1\endsubhead
Let $w\in W,E\in\Irr(W)$. We show:

(a) {\it $|\cb^F|\z_{w,E,q}$ is a polynomial in $q$ of degree $\le|w|+2\nu$;}

(b) {\it $|\cb^F|\z_{w,1,q}$ is a monic polynomial in $q$ of degree $|w|+2\nu$.}
\nl
From 1.5(b) and the definitions, for any $E'\in\Irr(W)$, 

(c) $c_{w,E',q}$ is a polynomial in $q$ of degree $\le|w|$;
\nl
moreover, we have 

$$\align&c_{w,1,q}=\sha(W)\i\sum_{\ce\in\Irr(\ch)}\sum_{z\in W}(\cf_\ce:R^1_z)_{G^F}\tr(T_w,\ce)\\&=
\sum_{\ce\in\Irr(\ch)}(\cf_\ce:1)_{G^F}\tr(T_w,\ce)=\tr(T_1,\ce_1)=q^{|w|}\tag d\endalign$$ 
where $\ce_1\in\Irr(\ch)$ is such that $T_y$ acts on it as $q^{|y|}$ for any $y\in W$.

For any $\g\in\un\cu$ and any $E'\in\Irr(W)$ we consider the sum
$$M_{E,E',\g}=\sum_{u\in\g^F}\tr(F,H^*(\cb_u)_{E'})\tr(F,H^*(\cb_u)_E).$$
(Notation of 1.5(a).) 

(e) This is a polynomial in $q$ of degree $\le2\dim\cb_u+\dim\g=2\nu$ (where $u\in\g$)
with the inequality being strict unless $E,E'$ appear in the top cohomology of $\cb_u,u\in\g$.
\nl
We have
$$\sha(\cb^F)\z_{w,q;E}=\sum_{E'\in\Irr(W)}c_{w,E',q}\sum_{\g\in\un\cu}M_{E,E',\g}.$$
(See 1.6(c).) Using this together with (e) and (c) we see that (a) holds. Now assume that $E=1$.
If $\g$ is not the regular unipotent class then $E$ does not appear in the top 
cohomology of $\cb_u,u\in\g$ hence by (e), $M_{E,E',\g}$ is a polynomial in $q$ of degree $<2\nu$.
If $\g$ is the regular unipotent class and $E'\ne1$ then $M_{E,E',\g}=0$.
If $\g$ is the regular unipotent class and $E'=1$ then $M_{E,E',\g}=|\g^F|$ is a monic polynomial in $q$
of degree $2\nu$. Combining this with (d) we see that (b) holds.

\mpb

Now let $C\in\un W$ and let $w\in C_{min}$. From (a),(b) we deduce

(f) If $i>|w|-(r-m(C))$ then $\Ps_i(C)=0$.

(g) If $i=|w|-(r-m(C))$ then $\Ps_i(C)\ne0$; more precisely, the multiplicity of $1$ in $\Ps_i(C)$ is $1$.
\nl
Note that $|w|-(r-m(C))$ is even. 
We conjecture that if $i\in\NN,i'\in\NN$ satisfy $i+i'=|w|-(r-m(C))$ then $\Ps_i(C)=\Psi_{i'}(C)$. In particular
the multiplicity of $1$ in $\Ps_0(C)$ is $1$. This is supported by the examples in \S3.

\subhead 2.2\endsubhead
For $\g\in\un\cu$ we set $d(\g)=\dim\cb_u$ where $u\in\g$.
According to Springer, for any $E\in\Irr(W)$ there is a unique $\g\in\un\cu$ such that
$H^{2d(\g)}(\cb_u)_E\ne0$ for some/any $u\in\g$; we set $\Xi(E)=\g$, $d'(E)=d(\g)$. The map 
$\Xi:\Irr(W)@>>>\un\cu$ is surjective. We shall need the following property of $\Xi$.

(a) {\it Let $\g\in\un\cu$, $u\in\g$ and let $E\in\Irr(W)$ be such that $H^i(\cb_u)_E\ne0$ for 
some $i$. Then $d'(E)=d(\Xi(E))\le d(\g)$. If in addition $i<2d(\g)$ then $d'(E)<d(\g)$.} 
\nl
This follows from 1.7(a).

\subhead 2.3\endsubhead
Let $\Ph:\un W@>>>\un\cu$ be the (surjective) map defined in \cite{L11}. By definition, if $C\in\un W$ and
$w\in C_{min}$, then  $\Ph(C)\cap(B^+wB^+)\ne\emp$; moreover, if $\g'\in\un\cu$ satisfies 
$\g'\cap(B^+wB^+)\ne\emp$ then $\Ph(C)$ is contained in the closure of $\g'$ in $\cu$.

Let $C\in\un W$ and let $\g=\Ph(C)$. Let $w\in C_{min}$. We show: 

(a) {\it If $E\in\Irr(W)$ appears in $\un\z_{C,q}$ then $d'(E)\le d(\g)$.}
\nl
From 1.3(b) we see that $H^i(\cb_u)_E\ne0$ for some $u\in\cu\cap(B^+wB^+)$. Now (a) follows from 2.2(a).

We show:

(b) {\it If $E\in\Irr(W)$ appears in $\un\z_{C,q}$ and $\Xi(E)\ne\g$ then $d'(E)<d(\g)$.}
\nl
From 1.3(b) we see that $H^i(\cb_u)_E\ne0$ for some $\g'\in\un\cu$ and some $u\in\g'\cap(B^+wB^+)$. 
If $i<2d(\g')$ then by 2.2(a) we have $d'(E)<d(\g')$. Using (a) we deduce $d(\g)<d(\g')$; but by 
the definition of $\Ph$, $\g$ is contained in the closure of $\g'$ so that $d(\g)\ge d(\g')$, a contradiction.
If $i=2d(\g')$ then $\Xi(E)=\g'$ hence $d'(E)=d(\g')$. Since $\g$ is contained in the closure of $\g'$ 
and $\g'\ne\g$ we have $d(g')<d(\g)$ hence $d'(E)<d(\g)$. This proves (b).

(c) {\it If $E_0\in\Irr(W)$ 
is the part of $H^{2d(\g)}(\cb_{u_0})$ which is fixed by the action of
the centralizer of $u_0$ ($u_0\in\g$) then $E_0$ appears in $\un\z_{C,q}$.}
\nl
We must show that $\sum_{u\in\cu^F\cap(B^+wB^+)}\tr(F,H^*(\cb_u)_{E_0})\ne0$.
Assume that 

(d) $\sum_{u\in\cu^F\cap(B^+wB^+)}\tr(F,H^*(\cb_u)_{E_0})=0$.
\nl
If $\g'\in\un\cu$ is such that $\g$ is contained in the closure of $\g'$ and $\g\ne\g'$ then
$H^i(\cb_u)_{E_0}=0$ for all $i$ hence from (d) we deduce 
$\sum_{u\in\g^F\cap(B^+wB^+)}\tr(F,H^*(\cb_u)_{E_0})=0$.
It follows that $\sha(u\in\g^F\cap(B^+wB^+))q^{d(\g)}=0$. This contradicts 
$\g^F\cap(B^+wB^+)\ne\emp$. This proves (c).

\subhead 2.4\endsubhead
Let $Z=\Ps(\un W)\sub\ZZ[\qq]\ot\car_W$. For $\x\in Z$ there is a unique $\g\in\un\cu$ such that the following 
holds: 

If $E\in\Irr(W)$ appears in $\x$ then $d'(E)\le d(\g)$. If $E\in\Irr(W)$ appears in $\x$ and $\Xi(E)\ne\g$ then 
$d'(E)<d(\g)$. If $E_0\in\Irr(W)$ is the part of $H^{2d(\g)}(\cb_{u_0})$ which is fixed by 
the action of the centralizer of $u_0$ ($u_0\in\g$), then $E_0$ appears in $\x$.

The existence of $\g$ follows from 2.3; the uniqueness is obvious. Thus $E\m\g$ is a well defined map
$\Th:Z@>>>\un\cu$. From 2.3 we have that $\Ph(C)=\Th(\Ps(C))$ for any $C\in\un W$. Since $\Ph$ is surjective, we
see that $\Th$ is surjective.

\subhead 2.5\endsubhead
Let $G_{ad}=G/Z_G$.
From 1.3(a) for any $u\in\cu^F$ and $E\in\Irr(W)$ we have
$$\tr(F,H^*(\cb_u)_E)=\sha(W)\i\sum_{y\in W}\tr(y,E)\tr(u,R^1_y).\tag a$$
For $y,y'$ in $W$ we have
$$\sum_{u\in\cu^F}\tr(u,R^1_y)\tr(u,R^1_{y'})=\sha(z\in W;zyz\i=y')\det(q-y,\r)\i\sha(G_{ad}^F),$$
see \cite{DL,Theorem 6.9}. Using this and (a) we obtain
$$\align&\sum_{u\in\cu^F}\tr(F,H^*(\cb_u)_{E'})\tr(F,H^*(\cb_u)_{E})\\&=
\sha(W)^{-2}\sum_{y\in W,z\in W}\tr(y,E')\tr(zyz\i,E)\det(q-y,\r)\i\sha(G_{ad}^F)\\&
=\sha(W)\i\sum_{y\in W}\tr(y,E\ot E')\det(q-y,\r)\i\sha(G_{ad}^F).\tag b\endalign$$
Using this and 1.6(c) we see that if $C\in\un W$, $w\in C_{min}$ then the coefficient of $E\in\Irr(W)$ in 
$\un\z_{C,q}$ is
$$\sum_{\ce\in\Irr(\ch),E'\in\Irr(W)}A_{C,\ce}A'_{\ce,E'}A''_{E',E}\tag c$$
where
$$A_{C,\ce}=(q-1)^{m(C)}\tr(T_w,\ce),$$
$$A'_{\ce,E'}=\sha(W)\i\sum_{z\in W}(\cf_\ce:R^1_z)_{G^F}\tr(z,E'),$$
$$A''_{E',E}=\sha(W)\i\sum_{y\in W}\tr(y,E\ot E')\det(q-y,\r)\i $$
Thus this coefficient is an entry of a product of three square matrices of size $\sha(\un W)=\sha(\Irr(W))$.
The matrices $(A_{C,\ce}),(A''_{E',E})$ are known to be invertible when regarded as matrices with entries in $\QQ$.
We see that 

(d) {\it $\{\un\z_{C,q};C\in\un W\}$ is a basis of $\QQ\ot\car_W$ if and only if the matrix
$(A'_{\ce,E'})$ (with entries independent of $q$) is invertible.}
\nl
The entries of the last matrix are some of the entries of the nonabelian Fourier
transform \cite{L84} for the various two-sided cells of $W$. This matrix is not necessarily
invertible. For example if $W$ is of type $B_n$ or $C_n$ this matrix is invertible if $n\le11$ but is not
invertible if $n=12$. Also if $W$ is of type $E_6$ this matrix is not invertible.

\subhead 2.6\endsubhead
Combining 1.5(a),1.6(c),2.5(b), we obtain the identity
$$\align&\sha(\cb^F)a_{w,y}=\sha(W)^{-2}\sum_{E,E',\ce\text{ in }\Irr(W)}\sum_{z,y'\text{ in }W}
\\&(\cf_\ce:R^1_z)\tr(z,E')\tr(T_w,\ce)
\tr(y',E)\tr(y',E')\det(q-y',\r)\i\sha(G_{ad}^F)\tr(y,E),\endalign$$
where $(\cf_\ce:?)$ is the multiplicity of
$\cf_\ce$ in the virtual representation $?$ of $G^F$, $\r$ is the reflection
representation of $W$ and $G_{ad}$ is the adjoint group of $G$.
We now replace $\sum_{E'\in\Irr(W)}\tr(z,E')\tr(y',E')=\sha(e\in W;eze\i=y')$ and we obtain
$$\align&\sha(\cb^F)a_{w,y}
=\sha(W)\i\sum_{E,\ce\text{ in }\Irr(W)}\sum_{z\in W}
\\&(\cf_\ce:R^1_z)\tr(T_w,\ce)\tr(z,E)\det(q-z,\r)\i\sha(G_{ad}^F)\tr(y,E).\endalign$$
We now replace $\sum_{E\in\Irr(W)}\tr(z,E)\tr(y,E)$ by $\sha(e'\in W;e'ze'{}\i=y)$ and we obtain
$$a_{w,y}=\sum_{\ce\in\Irr(W)}
(\cf_\ce:R^1_y)\tr(T_w,\ce)\det(q-y,\r)\i\sha(G_{ad}^F)
\sha(\cb^F)\i.$$
We now replace $\det(q-y,\r)\i\sha(G_{ad}^F)$ by $\dim(R^1_y)q^\nu
(-1)^{|y|}$ and we
obtain
\proclaim{Proposition 2.7} For any $y\in W$ we have
$$a_{w,y}=\sum_{\ce\in\Irr(W)}
(-1)^{|y|}(\cf_\ce:\dim(R^1_y)R^1_y)\tr(T_w,\ce)q^\nu\sha(\cb^F)\i.\tag a$$
\endproclaim

\subhead 2.8\endsubhead
In this subsection we diverge from the setup of 1.1; we assume instead that $W$ is a finite Coxeter group. 
Then all ingredients of 2.5(a) make sense. The (constant) matrix 
$(A'_{\ce,E'})$ makes sense as a matrix whose entries involve the appropriate generalization of the
nonabelian Fourier transform. 
Hence the map $\Ps$ can be defined in this generality (although the ring of coefficients $\ZZ$ may have to be
increased).

\head 3. Examples\endhead
\subhead 3.1\endsubhead
We return to the setup in 1.1.
We shall denote by $\sgn$ the sign representation of $W$, by $1$ the unit representation of $W$. Let $\r$ be 
as in 1.8. In this section we shall sometime write $w$ instead of $C$ when $w$ is an element of $C\in\un W$.

Assume first that $W$ is of type $A_1$. The elements of $\un W$ are represented by $1,s$
where $s$ is the simple reflection. The objects of $\Irr(W)$ are $1,\r=\sgn$. We have

$\Ps(1)=1+\sgn$,

$\Ps(s)=1$.

\subhead 3.2\endsubhead
We now assume that $W$ is of type $A_2$. The elements of $\un W$ are represented by $1,s_1,c$
where $s_1,s_2$ are the simple reflections and $c=s_1s_2$. The objects of $\Irr(W)$ are $1,\r,\sgn$. We have

$\Ps(1)=1+2\r+\sgn$,

$\Ps(s_1)=1+\r$,

$\Ps(c)=1$.

\subhead 3.3\endsubhead
We now assume that $W$ is of type $B_2$. The elements of $\un W$ are represented by $1,s_1,s_2,c,c^2$
where $s_1,s_2$ are the simple reflections and $c=s_1s_2$. The objects of $\Irr(W)$ are
$1,\r,\e',\e'',\sgn$ where $\e',e''$ are the one dimensional representations other than
$1,\sgn$. The following result was obtained by making use of \cite{SR}. We can arrange notation so that:

$\Ps(1)=1+2\r+\e'+\e''+\sgn$,

$\Ps(s_1)=1+\r+\e'$,

$\Ps(s_2)=1+\r+\e''$,

$\Ps(c)=1$,

$\Ps(c^2)=\qq^21+\qq\r+1$.

\subhead 3.4\endsubhead
We now assume that $W$ is of type $G_2$. The elements of $\un W$ are represented by $1,s_1,s_2,c,c^2,c^3$
where $s_1,s_2$ are the simple reflections and $c=s_1s_2$. The objects of $\Irr(W)$ are
$1,\r,\e',\e'',\r'=\r\ot\e'=\r\ot\e'',\sgn$ where $\e',e''$ are the one dimensional representations other than
$1,\sgn$. The following result was obtained by making use of \cite{CR}. We can arrange notation so that:

$\Ps(1)=1+2\r+2\r'+\e'+\e''+\sgn$,

$\Ps(s_1)=1+\r+\r'+\e'$,

$\Ps(s_2)=1+\r+\r'+\e''$,

$\Ps(c)=1$,

$\Ps(c^2)=\qq^21+\qq\r+1$,

$\Ps(c^3)=\qq^41+\qq^3\r+\qq^2(\r'+1)+\qq\r+1$.

\subhead 3.5\endsubhead
In the examples above, $\Ps_i(C)$ is always an actual representation of $W$. This is not so in higher rank.
We have calculated $\Ps(C)$ for $G$ of type $F_4$ on a computer using the formulas in 2.5, by the same method 
as in \cite{L11, 1.2}; we thank Gongqin Li for programming this in GAP.
If $C$ consists of the longest element in $W$ of type $F_4$ and $E$ is a one-dimensional representation of $W$
other than $1,\sgn$ then the coefficient of $E$ in $\Ph(C)$ is $-2\qq^{10}$; thus, $\Ph_{10}(C)$ is not an actual
representation of $W$. A similar thing happens in type $B_3$ and $C_3$.

\enddocument